\documentclass[12pt]{amsart}
\usepackage{amsthm, amsfonts, amssymb, amsmath}
\DeclareMathOperator{\res}{Res}

\DeclareMathOperator{\Span}{span}
\DeclareMathOperator{\reg}{reg}
\DeclareMathOperator{\gr}{gin_{revlex}}
\DeclareMathOperator{\gl}{gin_{lex}}

\DeclareMathOperator{\In}{in_{\tau}}
\DeclareMathOperator{\syl}{Syl}

\DeclareMathOperator{\sat}{sat}
\DeclareMathOperator{\xdeg}{deg_{x_0}}
\DeclareMathOperator{\Seg}{Seg}
\DeclareMathOperator{\Lex}{Lex}
\DeclareMathOperator{\gin}{gin}

\DeclareMathOperator{\lex}{lex}
\DeclareMathOperator{\revlex}{revlex}
\DeclareMathOperator{\ini}{in}
\DeclareMathOperator{\incoef}{incoef}
\DeclareMathOperator{\GL}{GL}
\DeclareMathOperator{\sgn}{sign}

\def\PP{\mathbb{P}}
\def\NN{\mathbb{N}}

\newtheorem{thm}{Theorem}[section]
\newtheorem{cor}[thm]{Corollary}
\newtheorem{lem}[thm]{Lemma}

\newtheorem{prop}[thm]{Proposition}
\theoremstyle{definition}
\newtheorem{ex}[thm]{Example}
\newtheorem{defn}[thm]{Definition}
\newtheorem{rmk}[thm]{Remark}

\begin{document}
\title{Generic initial ideals of points and curves}
\author{Aldo Conca}
\email{conca@dima.unige.it}
\address{Dipartimento di Matematica, Universit\'a di Genova, Via Dodecaneso 35,  I-16146 Genova, Italia}
\author{Jessica Sidman}
\email{jsidman@mtholyoke.edu}
\thanks{The second author was supported by an NSF postdoctoral fellowship during 2002-2003.}
\address{415a Clapp Lab, Department of Mathematics and Statistics, Mount Holyoke College, South Hadley,
MA 01075.}
\begin{abstract} Let $I$ be the defining ideal of a smooth complete intersection space
curve $C$ with defining equations of degrees $a$ and $b$. We use the partial elimination ideals
introduced by Mark Green to show that the lexicographic generic initial ideal of $I$  has
Castelnuovo-Mumford regularity $1+ab(a-1)(b-1)/2$  with the exception of the case
$a=b=2,$ where the regularity is   $4$.  Note that $ab(a-1)(b-1)/2$ is exactly   the number of singular
points of a general projection of $C$ to the plane.   Additionally, we show that for any term ordering
$\tau$, the generic initial ideal of a generic set of points in $\PP^r$ is a $\tau$-segment ideal. 
\end{abstract}
\maketitle

\section{Introduction}\label{sec: intro}

 Let $S = k[x_0, \ldots, x_r]$ where $k$ is an algebraically closed field
of characteristic zero and let  $\tau$ be   a term ordering on $S$.
Let $I \subset S$ be  a homogeneous ideal.
There is a monomial ideal canonically associated with 
$I$, its generic initial ideal with respect to $\tau$ denoted by
 $\gin_\tau(I)$ or simply  $\gin_{\tau} I.$ 
  In  this paper we study lexicographic 
generic initial ideals of curves and points 
via Green's partial  elimination ideals.  

For a  smooth    complete intersection curve $C$ in $\PP^3,$  we show that the complexity of 
its lexicographic generic initial ideal,   as measured by Castelnuovo-Mumford regularity, is
 governed by the geometry of a generic projection of $C$ to $\PP^2$.

\begin{thm}\label{thm: curves}  Let $C$ be a smooth    complete intersection of  
hypersurfaces of degrees
$a,b > 1$ in $\PP^3$.    The regularity of the  lexicographic generic initial ideal of $C$ is equal to 
\[ 
\left \{
\begin{array}{ll} 1+\frac{a(a-1)b(b-1)}{2} & \mbox{  if } (a,b)\neq (2,2)\\  \\ 4                       
& \mbox{  if } (a,b) = (2,2)
\end{array} 
\right. \]  
\end{thm}

Note that, apart for the exceptional case $a=b=2$,  the regularity of the lexicographic generic initial
ideal is $1+$ the number of nodes of the generic projection of $C$ to $\PP^2.$  The statement of Theorem \ref{thm: curves} generalizes Example 6.10 in \cite{green}
which treats the special case where $a = b = 3.$

Macaulay's  characterization of Hilbert functions, see for instance Theorem 4.2.10 in \cite{bruns-herzog},  implies that
any  ideal $J$ is generated in degrees bounded by the largest degree of a 
generator of the corresponding lex-segment $\Lex(J)$.  
Much more is true -- Bigatti \cite{bigatti}, Hullett  \cite{hulett} and Pardue \cite{pardue} showed the Betti numbers
of $J$ are 
bounded by those of $\Lex(J)$. 
Let $I$ be the ideal of $C$ in Theorem \ref{thm: curves}.
For such an ideal $I$ one can  compute the largest degree of a 
generator of $\Lex(I)$. This  has been  done, for instance, by  D. 
Bayer in his Ph.D.  thesis   (Proposition in II.10.4, \cite{bayer})    
and by Chardin and Moreno -Socias \cite{CM}, and it turns out to be   
$\frac{a(a-1)b(b-1)}{2}+ab$. 
So the lexicographic generic initial ideal in Theorem 
\ref{thm: curves} is not equal to the lex-segment ideal but nearly achieves the worst-case regularity for its
Hilbert  function.  Moreover, as Bermejo and Lejeune-Jalabert have shown in \cite{bermejo-lejeune}, the 
extremal bound can only be achieved if $C$ lies in a plane.

We also study the generic initial ideals of finite sets of points. Surprisingly, when $X$ is a 
set of generic points its generic initial ideal is an \emph{initial segment}.

\begin{thm}\label{thm:points}
 Let $I$ be the ideal of $s$ generic points  of $\PP^n.$ Then  $\gin_\tau I$ is equal to the
$\tau$-segment ideal $\Seg_\tau(I)$ for all term orders $\tau$. In particular,
$\gin_{\lex} I$ is a lex-segment ideal. 
\end{thm}

The genericity required in Theorem \ref{thm:points} is quite explicit: the conclusion holds for a set
$X$ of $s$ points  if there is a system of coordinates such that the defining ideal of
$X$ does not contain non-zero forms supported on $\leq s$ monomials.  A special
case of the result when $\tau = \revlex$ is proved by Marinari and Ramella in \cite{marinari-ramella}.

For an introduction to generic initial ideals see \S 15.9 in \cite{eisenbud}. Here we just recall:

\begin{thm}[Galligo, Bayer-Stillman]Given a homogeneous ideal $I$ and a term ordering $\tau$ on the  monomials of
$S,$ there exists a dense open subset  $U \subseteq \GL_{r+1}(k)$
 such that $\gin_{\tau} I := \In (g\cdot I)$ is constant over all $g \in U$
and $\gin_{\tau}I$ is Borel-fixed. 
\end{thm}

Recall also that, in characteristic $0$,   an ideal $J$ is Borel-fixed if it is monomial and satisfies:

\[ \mathrm{if \ }m \ \mathrm{is \ a \ monomial,  \ } x_im \in J \implies x_jm \in J, \ \forall j \leq i.\]

From this property  one easily shows that the regularity of a Borel-fixed ideal $J$ is the maximum degree of a 
minimal generators.  A minimal resolution of such an ideal was constructed by Eliahou and Kervaire in
\cite{eliahou-kervaire}.

As $\tau$ varies over all term orderings, both the regularity and the minimal number of generators of
$\gin_{\tau} I$ may vary greatly.  The generic initial ideals with respect to the reverse lexicographic 
(revlex) term ordering have the minimum level of complexity possible.

\begin{thm}[Bayer-Stillman \cite{bayer-stillman}] If $I$ is a homogeneous ideal of $S,$ and $J= \gr I$  then
\[\reg I = \reg  J=\mbox{ max degree of a minimal generator of }  J.\]
\end{thm}    

The paper is organized as follows.  In \S 2, we setup notation and review terminology.  We introduce
partial elimination ideals, their basic properties, and algorithms for their computation in \S 3.  We
focus on the case of complete intersection curves in \S 4 and on the case of points in \S 5.
\vspace{.5 cm}

\noindent {\bf Acknowledgements} We are deeply indebted to D. Speyer as much of this work originated
jointly with him.  We thank B. Sturmfels for his support and encouragement.  We would also like to thank
D. Eisenbud, R. Thomas, and M. Stillman for helpful discussions regarding partial elimination ideals and
their computation and A. Dickenstein, L. Bus\'e, and C. D'Andrea for many stimulating conversations
about subresultants.  We have benefitted substantially from experimentation with the computer algebra
packages \emph{Macaulay 2} \cite{grayson-stillman} and \emph{CoCoA}
\cite{cocoa}.  We would like to thank R. Lazarsfeld and anonymous referees
for helpful comments.  We are also grateful for the support of MSRI, and the second author thanks the Clare Boothe Luce Program.

\section{Notation and terminology}  Let $S = k[x_0, \ldots, x_r]$ where $k$ is an algebraically closed field
of characteristic zero. Denote by  $\mathfrak{m}$ the irrelevant maximal ideal of $S$. 
  For an element $\alpha = (\alpha_0, \ldots, \alpha_r) \in \NN^{r+1}$ we let
$x^{\alpha}$ denote $x_0^{\alpha_0}\cdots x_r^{\alpha_r}.$  In this section we briefly recall notions
related to term orderings and Castelnuovo-Mumford regularity.  For a comprehensive introduction to general
notions related to Gr\"obner bases see \cite{cox-little-oshea} and \cite{kreuzer-robbiano}.

\begin{defn}\label{defn: orders} We say that a total ordering $\tau$ on the monomials of $S$ is a
\emph{term ordering} if it is a well-ordering satisfying \[x^{\alpha} >_{\tau} x^{\beta} \ \Rightarrow \
x^{\gamma} \cdot x^{\alpha} >_{\tau} x^{\gamma} \cdot x^{\beta} \ \ \forall \alpha, \beta, \gamma \in
\NN^{r+1}.\]  \end{defn} 

A term ordering $\tau$ on $S$ allows us to assign to each non-zero element $f \in S$ an \emph{initial term}
$\ini_\tau(f)$ and to any ideal $I$ an initial ideal $\ini_\tau(I)$. 

In what follows we will work exclusively with homogneous ideals and  we will always require that  the
term ordering is degree compatible: $m>n$ if $\deg(m)>\deg(n)$.

The \emph{lexicographic} and (degree) \emph{reverse lexicographic} term orderings feature prominently in
the literature.  If $x^{\alpha}$ and $x^{\beta}$ are two monomials of the same degree, then $x^{\alpha}
>_{\lex} x^{\beta}$ if the left-most non-zero entry of $\alpha-\beta$ is positive and $x^{\alpha}
>_{\revlex} x^{\beta}$ if the right-most non-zero entry of $\alpha-\beta$ is negative.  

Although our primary motivation for studying partial elimination ideals is to understand lexicographic
initial ideals, partial elimination ideals also provide a mechanism for studying initial ideals with
respect to any \emph{elimination} order.

\begin{defn} An \emph{elimination order} for the first $t$  variables  of $S$ is a term order $\tau$
such that if $f$ is a polynomial whose initial term $\In(f)$ does not involve variables $x_0,
\ldots,x_{t-1},$  then $f$ itself does not involve  variables $x_0, \ldots,x_{t-1}.$ 
\end{defn}

As we shall see in Proposition \ref{prop: gb calculation of K_p}, one may use an elimination order for
the variable
$x_0$ to compute partial elimination ideals.  If $\tau$ is an elimination order for $x_0,$ then it
is equivalent to a $(1,r)$ \emph{product order}  which first sorts monomials by powers of $x_0$ and then
sorts the remaining variables by an arbitrary term ordering
$\tau_0.$

We will use the notion of \emph{Castelnuovo-Mumford regularity} as a rough measure of the complexity of
our computations.

\begin{defn}\label{defn: reg} Let $M$ be a finitely generated graded $S$-module, and let 
\[ 0 \to \oplus_j S(-a_{lj}) \to  \cdots \oplus_j S(-a_{1j}) \to \oplus_j S(-a_{0j}) 
\to M \to 0\] be a minimal graded free resolution of $M.$  We say that $M$
 is \emph{$d$-regular} if $a_{ij} \leq d+i$ for all $i$ and $j,$ and that  the \emph{regularity} of $M,$
denoted $\reg M,$ is the least $d$ such that
$M$ is $d$-regular.
\end{defn}

One may also formulate the definition of regularity in terms of vanishings of local cohomology with
respect to $\mathfrak{m}.$  The vanishing of the zero-th local cohomology group is related to the notion
of \emph{saturation} which plays an important role in the study of regularity. 

\begin{defn}\label{defn: sat} Let $I \subseteq S$ be a homogeneous ideal.  The \emph{saturation} of $I,$ 
denoted $I^{\sat}$ is defined to be $I:_S \mathfrak{m}^{\infty}.$  Note that
$I_d = I^{\sat}_d$ for all $d \gg 0.$  We say that $I$ is $d$-saturated if $I$ agrees with its
saturation in degrees $d$ and higher.  The minimum degree for which $I$ is $d$-saturated is the
\emph{saturation degree} (also the \emph{satiety index} in \cite{green}) of $I.$  
\end{defn}

\section{Partial elimination ideals}
\label{sec: partial elimination ideals}  Let $S = k[x_0, \ldots, x_r],$ and let $\overline{S} = k[x_1,
\ldots, x_r].$    Let $\tau$ be an arbitrary elimination order on
$S$ that eliminates the variable $x_0$ and hence induces a term order,  denoted by $\tau_0$, on
$\overline{S}.$  In this  section we set up the theory of partial elimination ideals over a  polynomial
ring in $r+1$ variables as introduced in \cite{green}.  Much of the material in \S \ref{sec: basics} and \S \ref{sec: computation}
 appears either explicitly or implicitly in 
\cite{green}, but we give proofs here both to keep the presentation self-contained and to present a more algebraic point of view.

We represent any non-zero polynomial $f$ in $S$ as  
$$f=f_0x_0^p+f_1x_0^{p-1}+\dots+f_p$$  with
$f_i\in \overline{S}$ and $f_0\neq 0$. The polynomial $f_0$ is called the initial coefficient of
$f$ with respect to $x_0$ and is denoted by $\incoef_{x_0}(f).$  The integer $p$ is called the
$x_0$-degree of
$f$ and is denoted by $\deg_{x_0}(f)$.  

\subsection{Definitions and basic facts}\label{sec: basics} In this section we define the partial
elimination ideals and describe their  basic algebraic and geometric properties.  We begin with the
definition:

\begin{defn}[Definition 6.1 in \cite{green}]  Let $I$ be a homogeneous ideal in $S.$  The \emph{$p$-th
partial elimination ideal} of $I$ is defined to be the ideal

\[K_p(I) := \{ \incoef(f)\   | \  f \in I \mbox{ and }  \deg_{x_0} f=p \}\cup \{0\}\]
in the polynomial ring $\overline{S} = k[x_1, \ldots, x_r].$
 \end{defn}

It is easy to see that if $I$ is homogeneous then $K_p(I)$ is also homogeneous. 

In Lemma \ref{partialfacts} we gather together some elementary algebraic facts about
 the partial elimination ideals.  We leave the proof to the reader.  The decomposition of $\In I$ given in part  (1) is one of the
motivations for the definition.

\begin{lem}\label{partialfacts} Let $I$ be a homogeneous ideal.
\begin{enumerate}
\item  $ \In  I = \sum_p  x_0^p \ini_{\tau_0} K_p(I).$

\item Taking $K_p$ commutes with taking initial ideals:  $K_p(\ini_\tau I)=\ini_{\tau_0} K_p(I)$

\item The partial elimination ideals are an ascending chain of ideals, i.e., $K_i(I) \subseteq K_{i+1}(I)$ for
all $i$.

\end{enumerate}
\end{lem}

One expects that if $I$ is in generic coordinates, then the   partial elimination ideals $K_p(I)$ are
already in generic coordinates.   Proposition \ref{prop: gin of partial} shows that this is indeed the
case.

\begin{prop}\label{prop: gin of partial} Let $I \subset S$ be a homogeneous ideal.   If $I$ is in
generic coordinates  then $\ini_{\tau_0} K_p(I) = \gin_{\tau_0} K_p(I).$  
\end{prop}

\begin{proof} Let $\GL_{r}(k)$ act on $\overline{S}$ in the usual way and extend this to an action on
$S$ in the trivial fashion by letting elements of $\GL_{r}(k)$ fix $x_0.$  

We know that the ideal $I$ determines a dense open subset $U
\subset \GL_{r+1}(k)$ with the property that $g \in U$ implies that $\In (g I) = \gin_{\tau}I.$  We show 
that for each $g \in U$ there is a dense open subet $U' \subset \GL_{r}(k)$ so that for all $h \in U'$
\begin{enumerate}
\item $\gin_{\tau_0}(K_p(g   I))= \ini_{\tau_0}(h K_p(g   I))$ 
\item $hg$ is again a generic change of coordinates for $I.$
\end{enumerate}
  Consider the space $\GL_r(k) \times \GL_{r+1}(k)$ with projection maps $\pi_1$ and $\pi_2$ onto the
first and second factors, respectively.   The map \[\phi: \GL_r(k) \times \GL_{r+1}(k) \to
\GL_{r+1}(k)\] given by $\phi(h,g) = hg$ is regular.    The
inverse image of $U$ under the map $\phi$ is a dense open subset of $\GL_r(k) \times \GL_{r+1}(k).$ For
each $g \in U$ the set $W := \pi_1 ( \pi_2^{-1}(g) \cap \phi^{-1}(U))$ is a dense open 
subset of $\GL_r(k).$  The element $g$ determines a dense open set $V \subset \GL_r(k)$ such that $h \in V$
satisfies (2), i.e., each $h \in V$ is a set of generic coordinates for $K_p(g I).$  Then any $h \in U' := W
\cap V$ has the property that $hg$ is a set of generic coordinates for $I.$

For $h$ and $g$ chosen as above, we have $h K_p(g I)=K_p(hgI).$  Thus, 
$\gin_{\tau_0}(K_p(g   I))=  \ini_{\tau_0} K_p(hg   I).$   By Lemma
 \ref{partialfacts} (2), $\ini_{\tau_0} K_p(I)= K_p(\In I),$ which implies that
$\gin_{\tau_0}(K_p(gI))=K_p \ini_{\tau}(hgI).$

Since $hg$ is again generic, $\gin_{\tau} I= 
\ini_{\tau} gI= \ini_{\tau} hgI$.  So we have $\gin_{\tau_0}(K_p(gI))=K_p (\ini_{\tau}gI)$. Using Lemma 
\ref{partialfacts} (2) again, we obtain $\gin_{\tau_0}(K_p(gI))=\ini_{\tau_0} K_p(gI)$ and this proves
the assertion. 
\end{proof}

The partial elimination ideals of an arbitrary homogeneous ideal $I$ can be  recovered in an easy way
from a Gr\"obner basis for $I.$  In practice one may want to take a $(1,r)$ product order with the
reverse lexicographic ordering on the last
$r$ variables in order to minimize computations.

\begin{prop}\label{prop: gb calculation of K_p}   Let $G$ be a Gr\"obner basis for $I$ with respect to an
elimination ordering $\tau$.  Then the set 
$$ G_p=\{ \incoef_{x_0}(g) \mid g \in G  \mbox{ and } \deg_{x_0}(g)\leq p\}$$ is a Gr\"obner basis for
$K_p(I)$.   
\end{prop}

\begin{proof} Note that if $g\in I$ and $\deg_{x_0}(g)=p$ then $\incoef_{x_0}(g)\in K_p(I)$ by
definition. By Lemma 
\ref{partialfacts} (3) we have that the elements of $G_p$ are in $K_p(I)$.   We will show that
their initial terms   generate $\ini_{\tau_0} K_p(I)$.   Suppose that $m$ is a monomial in the ideal 
$\ini_{\tau_0}  K_p(I)$. This implies that there exists $f\in I$ such that
$\ini_\tau(f)=mx_0^p$ and hence  there exists $g\in G$ such that $\ini_\tau(g)| \ini_\tau(f)$. Set
$h=\incoef_{x_0} (g)$.  It follows that $\deg_{x_0} g\leq p$, so that $h\in G_p$,  and $\ini_{\tau_0} h
|m$.  
\end{proof} 

By  part (3) of Lemma \ref{partialfacts} we know that the subscheme cut out  by the $p$-th partial
elimination ideal is contained in the subscheme defined by
 the $(p-1)$-st partial elimination ideal.  The following result gives the  precise relationship between
the partial elimination ideals and the geometry of
 the projection map from $\PP^r$ to $\PP^{r-1}$.

\begin{thm}[Proposition 6.2 in \cite{green}]\label{thm: green1} Let $Z$ be a reduced subscheme of
$\PP^r$ not containing $[1:0: \cdots : 0]$  and let $I = I(Z)$ be the homogeneous ideal of $Z.$  Let 
\[\pi: \PP^r \to \PP^{r-1}\] be the projection from the point 
$[1:0: \cdots :0].$ Set-theoretically, $K_p(I)$ is the ideal of 
\[ \{z \in \pi(Z) \mid | \pi^{-1}(z)| > p\},\]
where $|\pi^{-1}(z)|$ denotes the length of the scheme-theoretic fiber above $p.$
\end{thm}
\begin{proof} We prove the theorem by reducing to the affine case.  We begin by introducing some notation.  If $J \subseteq S$ is a homogeneous ideal, let $J_{(x_i)}$
denote its dehomogenization in $k[ \frac{x_0}{x_i}, \frac{x_1}{x_i}, \ldots, \frac{x_r}{x_i}].$

To show that $K_p(I)$ cuts out the $(p+1)-$fold points set-theoretically it
suffices to show that $K_p(I)_{(x_i)}$ cuts out  the $(p+1)-$fold points in each of the standard
affine open patches of 
$\PP^{r-1}$ for $i = 1, \ldots, r.$  If we consider the ideal $I_{(x_i)} \subseteq k[ \frac{x_0}{x_i}, \frac{x_1}{x_i}, \ldots, 
\frac{x_r}{x_i}]$ with the term ordering induced by $\tau$ in the natural way on the monomials in 
 $\frac{x_0}{x_i}, 
\frac{x_1}{x_i}, \ldots, \frac{x_r}{x_i},$ then by Lemma 4.8.3 in \cite{haiman},
 $K_p(I_{(x_i)})$ is set-theoretically the ideal of the $(p+1)-$fold points lying in this affine
patch.  

It remains for us to show that for any 
$i = 1,\ldots, r,$ \[K_p(I)_{(x_i)} = K_p(I_{(x_i)}).\]

It is clear that 
$K_p(I)_{(x_i)} \subseteq  K_p(I_{(x_i)}).$  For the opposite inclusion, note
that $\frac{x_0}{x_i}$ appears in the dehomogenization of a monomial $m$ precisely
as many times as $x_0$ appears in $m,$ and apply the definitions.
\end{proof}

In the situation of Theorem \ref{thm: green1}, we can see that $K_0(I)$ is in fact radical.  The ideal
$K_0(I)$ is just equal to $I \cap \overline{S}$. On the other hand, the higher 
 $K_p(I)$ need not be radical even if $I$ is a prime complete intersection of codimension $2$  in
generic coordinates; see Example \ref{nonsmooth}.

\subsection{Partial elimination ideals for codimension $2$ complete intersection. }
\label{sec: computation}
\setcounter{MaxMatrixCols}{20}

 Let \[f = x_0^a + f_1x_0^{a-1}+\cdots +f_{a-1}x_0 +f_a\] and \[g = x_0^b + g_1x_0^{b-1}+ \cdots +
g_{b-1}x_0 +g_b\] where $f_1, \ldots, f_a$ and $g_1, \ldots, g_a$ are indeterminates.  

We wish to describe  the partial elimination ideals of the ideal $I_{a,b}$ generated by $f$ and $g$ in
$S = k[x_0, \ldots, x_r]$ after specializing the $f_i$ and the $g_i$ to homogeneous elements of
$\overline{S} = k[x_1, \ldots, x_r]$ of degree $i.$

As we saw in \S \ref{sec: basics}, the partial elimination ideals of an arbitrary  homogeneous ideal $I$
can be recovered from a Gr\"obner basis for
 $I$ and, vice versa, give information on that  Gr\"obner basis.  In this section we discuss a result of
Eisenbud and Green showing that  $K_p(I_{a,b})$ is generated by the minors of a truncation of the
Sylvester matrix
 as long as  the forms $f$ and $g$ are generic enough.   Both Theorem \ref{greensyl} and Lemma \ref{rank} are well-known to experts, but we give proofs for
completeness.

\begin{thm} [Proposition 6.8 (3) in \cite{green}]
\label{greensyl}
 Assume that the $f_i$ and the $g_j$ are independent indeterminates and that $p< a\leq b.$  Let \[R =
k[f_1, \ldots, f_a, g_1, \ldots, g_b, x_0],\] where $k$ is an arbitrary field.  Define
$\syl_p(f,g)$ to be the matrix consisting of the first
$a+b-p$ rows of the  Sylvester matrix of $f$ and $g$, i.e. 

\[
\syl_p(f,g) = \begin{pmatrix} 1    & 0       &        & 0        & 1    & 0       &        & 0\\ f_1   
& 1    &        & 0        & g_1    & 1    &        & 0  \\
\vdots & \vdots  & \ddots & \vdots   & \vdots & \vdots  & \ddots & \vdots\\ f_a    & f_{a-1} & \ddots &
1      &        &         & \ddots & 1 \\ 0      & f_a     & \ddots & \vdots   & 0      & g_b     &
\ddots & \vdots \\ 
\vdots & \vdots  & \ddots & f_{a-p-1}  & \vdots & \vdots  & \ddots & g_{b-p-1} \\ 0      & 0      
&        & f_{a-p}      & 0      & 0       &        & g_{b-p} 
\end{pmatrix}\]

Then the ideal $K_p(f,g) \subset R$ is generated by the maximal minors of the matrix $\syl_p(f,g)$. 
\end{thm}

\begin{proof} Let $R_{\leq t}$  denote the vector space of polynomials in $R$ with $\xdeg \leq t.$  To compute $K_p(f,g)$ we want to find all $A,B \in R$ such that 
\begin{equation}Af+Bg = c_0x_0^p +c_1x_0^{p-1} + \cdots + c_{p-1}x_0 + c_p.
\label{eq: p} \end{equation} with 
$c_i \in k[f_1, \ldots, f_a, g_1, \ldots, g_b].$  Note that it suffices to find all $A,B$  satisfying the equation (\ref{eq: p}) where
$\xdeg A \leq b-1$ and $\xdeg B \leq  a-1.$  

The matrix $\syl_p(f,g)$ gives a linear map 
\[R_{\leq b-1} \oplus R_{\leq a-1} \to R_{\leq a+b-1}/( 1, x_0, \ldots, x_0^{p-1}).\] The kernel of
$\syl_{p+1}(f,g)$ consists of the set of all $(A,B) \in R_{\leq b-1} \oplus R_{\leq a-1}$ that satisfy
equation  (\ref{eq: p}).  The image of
$\ker \syl_{p+1}(f,g)$ under $\syl_p(f,g)$ is exactly the set
\[ \{ c_0x_0^p \ | \ c_0x_0^p +c_1x_0^{p-1} + \cdots + c_{p-1}x_0 + c_p \in (f,g)\}.\]

We will show that the maximal minors of $\syl_p(f,g)$ generate the image of $\ker \syl_{p+1}(f,g)$ under
the  map $\syl_p(f,g)$ as long as $\syl_{p+1}(f,g)$ drops rank in  the expected codimension.  The proof
that $\syl_{p+1}(f,g)$ does indeed drop rank in codimension $p+2$ will be given in Lemma
\ref{rank}.

If $\syl_{p+1}(f,g)$ drops rank in the expected codimension, then since
 $R$ is Cohen-Macaulay we conclude that the Buchsbaum-Rim complex  resolves the cokernel of
$\syl_{p+1}(f,g).$  (See Eisenbud \cite{eisenbud} A2.6 for details.)  

Using the Buchsbaum-Rim complex we can give explicit formulas for elements of
$\ker \syl_{p+1}(f,g)$ indexed by $T \subseteq \{1, \ldots, a+b \}$ with $|T| = a+b-p.$   Define
$\syl_{p+1}(f,g)^T$ to be the $(a+b-p-1)\times(a+b-p)$ matrix consisting
 of all of the columns of $\syl_{p+1}(f,g)$ indexed by elements of $T.$  Define $W_T$ to be the vector
of length $a+b$ whose $i$-th entry is 0 if 
$i \notin T$ and $\sgn(i) \det \syl_{p+1}(f,g)^{T - \{ i\}}$ if $i \in T$ where
$\sgn(i) = 1$ if the number of elements of $T$ less than $i$ is even, and -1 if the number of elements of $T$ less than $i$ is odd.  The  Buchsbaum-Rim complex is a
resolution precisely when the vectors $W_T$ generate the kernel of $\syl_{p+1}(f,g).$

Finally, we apply $\syl_p(f,g)$ to the elements $W_T$ constructed above.  The  dot product of $W_T$ with
each of the first $a+b-p-1$ rows of $\syl_p(f,g)$ is  zero since $W_T$ is in the kernel of
$\syl_{p+1}(f,g).$  The dot product of 
$W_T$ with the last row is just the expansion of the maximal minor of 
$\syl_p(f,g)$ corresponding to the columns indexed by $T$ by this final row.  Therefore, 
\[\syl_p(f,g) \cdot W_T = \det \syl_p(f,g)^Tx_0^p.\]
\end{proof}

\begin{lem}\label{rank} If  $f_i$ and $g_j$ are independent indeterminates and $p< a\leq b$ the matrix 
$\syl_p(f,g)$ drops rank in the expected codimension  $p+1$.
\end{lem}

\begin{proof}  We will show that the set where $\syl_p(f,g)$ fails to have maximal rank, that is, where
$\dim_k \ker \syl_p(f,g) \ge p+1,$ has codimension $p+1$ in the space of all $f$ and $g$ where the $f_i$
and $g_i$ take values in $k.$   The result follows if we can show that for any specialization of the
indeterminates $f_i$ and $g_j$ to values in $k,$
$\dim_k \ker \syl_p(f,g) \ge p+1$ if and only if $f$ and $g$ have a common  factor of degree $p+1.$

It is clear that if $f$ and $g$ have a common factor of degree $p+1$ then 
$\dim_k \ker \syl_p(f,g) \ge p+1,$ since we can use the $(p+1)$ common factors
to construct $(p+1)$ syzygies on $f$ and $g$ with distinct degrees.  

To prove the other direction, we will use induction on $p.$  Suppose that $p = 0.$  Then 
\[\dim_k \ker \syl_0(f,g) > 0\] if and only if $\res(f,g) = 0.$  It is well-known (see \cite{cox-little-oshea}) that $\res(f,g)=0$ if
and only if $f$ and $g$ have a common factor of degree at least one.

We treat the case where $p>0.$ Our assumption implies that we can find $p+1$ linearly independent elements $(A_0, B_0), \ldots, (A_p, B_p)$ of the kernel of $\syl_p (f,g).$  Since \[A_0f+B_0g, \ldots, A_pf+B_pg \in \Span (1, \ldots, x_0^{p-1}),\] there
is a nontrivial linear relation $\sum \lambda_i(A_if +B_ig) = 0.$  Hence, $f$ and 
$g$ must have a common factor so that $f = (x-\alpha)f'$ and $g = (x-\alpha)g'.$
By induction, we will be done  if we can show that the
dimension of 
\[\{(A,B) \in k[x_0]_{\leq b-2}\oplus k[x_0]_{\leq a-2} \mid Af'+Bg' \in \Span (1, x_0, \ldots, x_0^{p-2})\}\] 
is
$\ge p.$  But, we can assume (after reordering and cancelling leading terms) 
that for $i\ge 1,$ $\deg A_i \leq b-2$ and $\deg B_i \leq a-2.$  Consequently,
for $i \ge 1,$
\[ A_if'+B_ig' \in \Span (1, x_0, \ldots, x_0^{p-2}).\] 

\end{proof} Note that the first $b$ rows of $\syl_p(f,g)$ contain constants.  Recall the  following fact:

\begin{lem}[pg. 10 in \cite{bruns-vetter}]\label{lem: minors} Suppose that $M = (m_{i,j})$ is a $p
\times q$ matrix with entries in a  commutative ring.  If $m_{p,q}$ is a unit, then the ideal generated
by the  maximal minors of $M$ is the same as the ideal generated by the maximal minors  of the $(p-1)
\times (q-1)$ matrix $N$ with entries 
\[n_{i,j} = m_{i,j} - m_{p,j}m_{i,q}m_{p,q}^{-1} \qquad  1\leq i\leq p-1, \ \  1\leq j\leq q-1 .\]
\end{lem}

We have the following:  

\begin{cor} [See the remark following Proposition 6.9 in 
\cite{green}.]
\label{smaller}  

Let $a\leq b$ and assume that the $f_i$ and $g_i$ are sufficently general  homogeneous polynomials of
degree $i$ in variables $x_1,\dots, x_r$.   Assume also that $p<a$.      
\begin{itemize}
\item[(1)] The ideal of maximal minors of $\syl_p(f,g)$ is always contained in $K_p(f,g)$. It has
the expected codimension, $p+1,$ if $p\leq r-1$.
\item[(2)] Assume $p\leq r-2$. Then we have: 
      \begin{itemize}  
\item[(a)]  $K_p(f,g)$ is equal to the ideal of maximal minors of $\syl_p(f,g)$.    
\item[(b)]  $K_p(f,g)$ is   also the ideal generated by the maximal minors of a matrix  of  size  $(a-p)
\times a$ whose $(i,j)-$th  entry is either $0$ or has degree $b+i-j$. 
\item[(c)] $\reg K_p(f,g)=ab+{a-p+1 \choose 2}- {a+1 \choose 2} + p (a-p-1)$. 
\end{itemize}  
\end{itemize}  
\end{cor}  

\begin{proof}  Let $R = k[f_1, \ldots, f_a, g_1, \ldots, g_b, x_0]$ where the $f_i$ and $g_j$ are
indeterminates as in Theorem \ref{greensyl}.  Generators for $K_p(f,g)$ as an ideal in $R$ also generate
the $p$-th partial elimination ideal of the ideal generated by $f$ and $g$ in the ring $R \otimes_k
k[x_1, \ldots, x_r],$ which we will denote by $K_p(f,g) \otimes k[x_1, \ldots, x_r].$  An elementary
argument shows that if $p+1 \leq r,$ then for sufficiently general forms $f_i, g_j \in k[x_1, \ldots,
x_r],$ the specialization of the matrix $\syl_p(f,g)$ still drops rank in the expected codimension.  

Thus, (1) and (2,a)  follows from the proof of Theorem \ref{greensyl} and from Lemma \ref{rank}. Part
(b) of (2) follows from (2, a) and from iterated use of Lemma \ref{lem: minors}. Finally (2, c) follows
from (2, b) and from Lemma \ref{EN-reg}. 
\end{proof} 

\begin{rmk} The above corollary is sharp in the sense that, in general,  $K_{r-1}(f,g)$ strictly
contains  the ideal of maximal minors of $\syl_{r-1}(f,g)$. For instance, one can check with CoCoA that 
this happens if $r=3$ and $a=b=4$.   
\end{rmk} 

\begin{lem}\label{EN-reg}  Let $X=(h_{ij})$  be an $m\times n$ matrix of forms with $m\leq n$. Assume
$a_1,\dots, a_m$ and
$b_1,\dots, b_n$ are integers such that $\deg(h_{ij})=a_i+b_j>0$ whenever $h_{ij}\neq 0$. Assume    that
the ideal $I_m$ of maximal minors of $X$ has the expected codimension $n-m+1$. Then 
$$\reg I_m=\sum_i a_i +\sum_j b_j +(\max(a_i)-1)(n-m)$$
\end{lem} 
\begin{proof} The Eagon-Northcott complex gives a resolution of $I_m$ which is minimal since the entries
of the matrices in the resolution are (up to sign) the entries of $X$ and $0$. Keeping track of the
shifts one obtains the formula above.  The same formula can be derived  from the result 
\cite[Cor.1.5]{bruns-herzog1} of Bruns and Herzog. Another formula for the regularity appears   in
\cite{BudCasGor}. 
\end{proof} 

In particular we have:

\begin{cor}  
\label{finalcor}  
Let $I$ be the ideal of a smooth   complete intersection $C$ in $\PP^3$ defined by two forms
$f$ and $g$ of degrees $a,b>1.$  Assume that $I$ is in generic coordinates.  We have: 
\begin{itemize}
\item[(a)]  $K_1(f,g)$ is equal to the ideal of maximal minors of $\syl_1(f,g)$ and has codimension
$2$ in $k[x_1,x_2,x_3].$    
\item[(b)]  $K_2(f,g)$  contains the ideal of maximal minors of $\syl_2(f,g)$ and both ideals have codimension
$3$ in $k[x_1,x_2,x_3].$    
\end{itemize}  
\end{cor}  

\begin{proof}  We will use a geometric argument to show that if $f$ and $g$ are in sufficiently general coordinates, then $\syl_1(f,g)$ has
codimension $2$ and $\syl_2(f,g)$ has codimension $3$ in $k[x_1,x_2,x_3]$.  Since these codimensions are the
expected values for those determinantal ideals, the conclusion will follow by Corollary \ref{smaller}. 

Recall the classical fact that a generic projection of a smooth irreducible curve in $\PP^3$ has only nodes as singularities.  (See Theorem IV.3.10 in
\cite{hartshorne}.) It follows that after a generic change of coordinates, the image of the projection from the point $[1:0:0:0]$ will have only nodes as singularities.  As a consequence, we see that for each point $q \in \PP^2,$ the fiber of the projection of the curve $C$ will contain at most two points, and the set of $q$ with $\pi^{-1}(q) = 2$ is finite.  In other words, 
$\deg \gcd( f(x_0, q), g(x_0,q)) \leq 2$  and equality holds for only finitely many $q.$  
From the proof of Lemma \ref{rank}, we can see $\syl_p(f, g)$ drops rank at $q$ if and only if $f(x_0, q)$ and  $g(x_0,q)$ have a common factor of degree $\ge p+1.$  Therefore, we see that $\syl_1$(f,g) drops rank at a finite set of points and and $\syl_2(f, g)$ does not drop rank at any point in $\PP^2.$
\end{proof}

\section{The lexicographic gin of a  complete intersection curve in $\PP^3$}

 Let $I_{a,b}$ be a codimension $2$ complete intersection ideal in the polynomial ring $S = k[x_0,x_1,x_2,x_3]$ defined by two forms
of degrees $a,b>1.$  Let $C = V(I_{a,b})$ be the curve  in
$\PP^3$ defined by $I_{a,b}$.  We will assume that   $C$ is smooth  and in generic
coordinates. In other words, we assume that $I_{a,b}$ is prime, that the singular locus of $S/I_{a,b}$ consists solely of the homogeneus maximal ideal and that $I_{a,b}$ is in
generic coordinates.    We have that
$C$ has degree
$ab$ and genus
$ab(a+b-4)/2 +1.$  From Theorem
\ref{thm:  green1} we know that $K_0(I_{a,b})$ is the radical ideal of the  projection  $\pi: C
\to \PP^2$ from the point $[1:0:0:0]$. Since $C$ is in generic coordinates by assumption, the
projection $\pi$ is generic. Proposition
\ref{prop: res} describes additional numerical data associated to $\pi(C).$

\begin{prop}\label{prop: res}
The ideal $K_0(I_{a,b})$ is generated by a single polynomial of degree $ab$. It cuts out  a degree $ab$ curve with
$ab(a-1)(b-1)/2$ nodes.
\end{prop}

\begin{proof} 
We already know that $K_0(I_{a,b})$ is the radical ideal of $\pi(C)$ which has degree $ab.$  So it remains to show that $\pi(C)$ has $ab(a-1)(b-1)/2$ nodes.

Since a general projection of any space curve has only nodes as singularities,  we have  that $\pi (C)$ is a plane curve with  only nodes as singularities.  Since
$C$ is the normalization of $\pi(C)$ and $C$ has genus
$ab(a+b-4)/2+1,$ $\pi(C)$ has $$\frac{(ab-1)(ab-2)}{2} - \left( \frac{ab(a+b-4)}{2}+1
\right)=\frac{a(a-1)b(b-1)}{2}$$ nodes (see Remark 3.11.1 in \cite{hartshorne}). 
\end{proof}

Already, we can begin to describe the generators of $\gl I_{a,b}:$
\begin{cor}
The ideal $\gl I_{a,b}$ contains $x_1^{ab}$ and this is the only generator that is not divisible by $x_0.$
\end{cor}

\begin{proof} The generators of $\gl I_{a,b}$ are elements of $x_0^p \gl K_p(I_{a,b})$ for various $p$.  
So  clearly, the generators of $\gl K_0(I_{a,b})$ are the only generators of 
$\gl I_{a,b}$ not containing a factor of $x_0.$  But $K_0(I_{a,b})$ is principal, generated by a form of
degree $ab$ in generic coordinates.   The leading term of such a form is
$x_1^{ab}.$
\end{proof}

We are ready to prove the main result of the paper:

\begin{proof} [Theorem \ref{thm: curves}]  Set $I=I_{a,b}$, $K_p=K_p(I)$.  By virtue of Lemma
\ref{partialfacts} and since $x_0^a\in
\gin_{\lex} I$  we have 
 $$\gin_{\lex} I=\sum_{p=0}^a  x_0^p  \gin_{\lex} K_p$$
From Proposition \ref{prop: res} we know that   $\gin_{\lex} K_0=(x_1^{ab})$.   The proof consists of three steps.  First,
 we compute the  regularity of $\gin_{\lex} K_1$ explicitly.  Then we show that the regularity of  $\gin_{\lex} K_p -p \leq 1 +
\reg \gin_{\lex} K_1$ for 
$2\leq p\leq a-1$.   Finally, we will show that $\gin_{\lex} I$ actually requires a generator of degree
$\frac{1}{2}a(a-1)b(b-1)+1,$ which will complete the proof.

By Corollary \ref{finalcor} we have that $K_1$ is the
ideal of maximal minors of a  matrix of size $(a-1)\times a$ whose $ij$ entry has degree $b+i-j$.  The
resolution of $K_1$ is  given by  the Hilbert-Burch complex. It is then easy to determine the  degree of
$K_1$ from the numerical data of the resolution.  We obtain that  $K_1$ is unmixed and of degree
$\frac{1}{2}a(a-1)b(b-1)$.  We also know that the radical of $K_1$ is the ideal of definition of  
$\frac{1}{2}a(a-1)b(b-1)$ points. It follows that $K_1$ itself is the radical ideal defining
$\frac{1}{2}a(a-1)b(b-1)$ points.  We can conclude  from  Corollary \ref{ginlexpoints}  that $\reg
\gin_{\lex}(K_1)=\frac{1}{2}a(a-1)b(b-1)$. 

We now prove that for $p >1,$ the degrees of the generators of $x_0^p \gin_{\lex} K_p$ are bounded above
by  $1+\reg \gin_{\lex} K_1$, that is, by  $1+\frac{1}{2}a(a-1)b(b-1)$. This will imply that $\reg
\gin_{\lex} I$ is $\max( ab, 1+\frac{1}{2}a(a-1)b(b-1)).$

 From Corollary \ref{finalcor} (2) we have that the ideal, say
$J$, of the maximal minors of
$\syl_2(f,g)$  is Artinian (i.e. $K[x_1,x_2,x_3]/J$ is Artinian) and is contained in $K_2$ and that  $J$
is contained in $K_p$ for $p>1$. The regularity of an Artinian ideal $D$  is given by  the smallest $k$
such that the $k$-th power of the maximal ideal is contained in the ideal $D$ and hence does not change
when passing to the initial ideal.   It follows that  $\reg \gin_{\lex} K_p \leq \reg J$ for every
$p>1$. So the generators of 
$x_0^p \gin_{\lex} K_p$ are  in degrees $\leq p+\reg J$.  Taking into consideration that $K_a=(1)$,  it
is enough to show that $\reg J\leq 
\frac{1}{2}a(a-1)b(b-1)+1-p$ for all
$p=2,\dots,a-1$. So we may assume
$a>2$ and we have to show that $\reg J\leq  \frac{1}{2}a(a-1)b(b-1)+2-a$. To compute the regularity of
$J$ we first use Lemma \ref{lem: minors} to get rid of the units in the matrix defining
$J$ and then we use Lemma \ref{EN-reg}.  We get $\reg(J)=ab+{a-1 \choose 2}- {a+1 \choose 2} + 2 (a-3)$.

So it remains to show that 

$$ ab+{a-1 \choose 2}- {a+1 \choose 2} + 2 (a-3) \leq \frac{1}{2}a(a-1)b(b-1)+2-a$$

that is 

$$1/2a^2b^2 - 1/2a^2b - 1/2ab^2 - 1/2ab - a + 7\geq  0$$

for  all $3\leq a\leq b$. This is a simple calculus exercise. 

To finish the proof, we will show that if $m$ is a minimal generator of $\gin_{\lex}K_1,$ of degree 
$\frac{1}{2}a(a-1)b(b-1),$ then $x_0m$ is a minimal generator of $\gin_{\lex} I.$  If $x_0m$ is not
a minimal generator of $\gin_{\lex} I,$ then it must be divisible by some monomial $n$ that is a minimal
generator of $\gin_{\lex} I.$  This implies that $n \mid x_0m$ and that $n$ must be in $\gin_{\lex} K_0.$
However, this means that $n \mid m$ and $n \in \gin_{\lex} K_1$ since $K_0 \subseteq K_1.$  This contradicts
our choice of $m$ as a minimal generator.  We conclude that $x_0m$ must be a minimal generator of $\gin_{\lex} I.$
\end{proof}

\begin{ex} 
\label{nonsmooth}
One can check (using CoCoA, for instance) that
$I=(x^3-yz^2,  y^3-z^2t)$  defines an irreducible complete intersection curve $C$ with just one singular
point and that  $K_1(gI)$  with
$g$ a generic change of coordinates  is not radical. Indeed, $K_1(gI)$  has degree $18$ and it defines
only $11$ points, namely the $11$ singular points of the generic projection of $C$ to $\PP^2$. In this
case, the regularity of $\gin_{\lex}(I)$ is $16$ and not $19$ as in the smooth case. 
\end{ex}  

\setcounter{MaxMatrixCols}{10}

\section{The regularity of gins of points }

Set $S = k[x_0, \ldots, x_r]$.  We start with the following well-known lemma:

\begin{lem}\label{1-dim} Let $I$ be a homogeneous ideal of $S$ such that $S/I$ has Krull  dimension $1$
and
$\deg(S/I)=e$. Set  $c=\min \{ j  |  \dim [S/I]_i=e \mbox{ for all }  i\geq j\}$. Then  $\reg(I)\leq
\max\{e , c\}$.
\end{lem}

\begin{proof}  Let $J$ be the saturation of $I$. Then $S/J$ is a 
$1$-dimensional CM (Cohen-Macaulay) algebra. It is well-known and easy to see that  $\reg(J)\leq 
\deg(S/J)=e$  and
$\dim [S/J]_i=e$ for all $i\geq e-1$.  Let $p$ denote the saturation degree (satiety index)  of $I$,
i.e. the least $j$ such that $I_i=J_i$ for all $i\geq j$.  From the  characterization of regularity  in
terms of local cohomology it follows immediately that 
$\reg(I)=\max\{\reg(J), p\}$. To conclude, it is enough to show that $p\leq \max\{e , c\}$. If $p>e$ 
then $I_i=J_i$ for all $i\geq p$  and
$I_{p-1}\subsetneq  J_{p-1}$. Thus, $\dim [S/I]_i=e$ for all 
$i\geq p$ and $\dim [S/I]_{p-1}>e$. Hence $p=c$ and we are done.
\end{proof}

\begin{cor}\label{1-dimCM} Let $I$ be a homogeneous ideal of $S$ such that $S/I$ has Krull  dimension
$1$. Assume that the Hilbert function of $I$ is equal to the Hilbert  function of a
$1$-dimensional CM ideal. (e.g. $I$ is an initial ideal of a 
$1$-dimensional CM ideal).   Then  $\reg(I)\leq  \deg(S/I)$.
\end{cor}

\begin{proof} This follows from Lemma \ref{1-dim} since  the assumption  implies that $\dim
[S/I]_i=\deg(S/I)$ for all $i\geq \deg(S/I)-1.$
\end{proof}

\begin{cor}\label{ginlexpoints} Let $I$ be the ideal of a set $X$ of  $s$ points of $\PP^r$. Then $$\reg
\gin_{\lex} I=s.$$
\end{cor}

\begin{proof} By Corollary \ref{1-dimCM} we have $\reg \gin_{\lex} I\leq s$.  A general projection of $X$ to
$\PP^1$ will give $s$ distinct points.  This implies that $x_{r-1}^s$ is in $\gin_{\lex} I$. Since we work
with the lex order, $x_{r-1}^s$ is a minimal generator of $\gin_{\lex} I$. 
\end{proof} 
 
We want to show now that for a set of generic points the gin lex and  indeed any gin has a very special
form: it is a segment ideal.  Consider the  polynomial ring
$S=k[x_0,\dots,x_r]$   equipped with a term order $\tau$. Assume that 
$x_0>_\tau x_1>_\tau \dots >_\tau x_r$.

\begin{defn}\label{defn:  seg} A vector space $V$ of forms of degree $d$ is said to be a $\tau$-segment
if it is generated  by monomials and for every monomial  $m$ in $V$ and every monomial $n$ of degree $d$
with $n>_\tau m$ one  has $n\in V$. 
\end{defn}

Given a non-negative integer $u\leq { r+d \choose r} $ there exists exactly  one
$\tau$-segment of forms of degree $d$ and of dimension $u$: it is the  space generated by the $u$ 
largest monomials of degree $d$ with respect to 
$\tau$ and it will be denoted by $\Seg_\tau(d,u)$. Given a homogeneous ideal
$I$ for every $d$ we  consider the $\tau$-segment $\Seg_\tau(d, \dim  I_d)$ and define

$$\Seg_\tau(I)=\oplus_d  \Seg_\tau(d, \dim I_d).$$

By the very definition, $\Seg_\tau(I)$ is a graded monomial vector  space and simple examples show that 
$\Seg_\tau(I)$ is not an ideal in general.  But  there are important exceptions:  Macaulay's numerical
characterization of   Hilbert functions
\cite[Thm.4.2.10]{bruns-herzog} can be rephrased by saying that for every homogeneous  ideal $I$ the
space 
$\Seg_{\lex}(I)$ is an ideal.     In the following lemma we
collect a few simple facts about segments that will be used in the proof of that result.  

\begin{lem}\label{segments}  Let $\tau$ be a term order and let $V\subset S_a$ be a $\tau$-segment 
with  $\dim S_a/V \leq a$. Then  $S_1V$ is a $\tau$-segment  with $\dim S_{a+1}/VS_1=\dim S_a/V$. 
\end{lem} 

\begin{proof} First observe that since $x_{r-1}^a> x_{r-1}^{a-j}x_r^j$ for $j=1,\dots, a$ we have that
$x_{r-1}^a\in V$ and hence $(x_0,\dots, x_{r-1})^a \subseteq V$.  To prove that $VS_1$ is a
$\tau$-segment assume that $n$ is a monomial of degree $a+1$ such that  $x_im<n$ with $m$ in $V$; we
have to show that
$n\in VS_1$. Let
$k$ be the largest  index such that $x_k$  divides  $n$, so that $n=x_kn_1$. If $k\geq i$ then
$x_in_1\geq x_k n_1>x_i m$. It follows that  $n_1> m$ and hence $n_1\in V$ so that $n\in VS_1$. If,
instead,  $k<i$ then
$n\in (x_0,\dots,x_{r-1})^{a+1}$ which is contained in $VS_1$ since we have seen already that 
$(x_0,\dots,x_{r-1})^{a}$ is contained in $V$. 

To conclude, it is enough to show that the map $\phi$
induced by multiplication by $x_r$ is  an isomorphism from  $ S_a/V$ to $S_{a+1}/VS_1$. We show first
that $\phi$  is injective. If $m$ is a monomial in $S_a\setminus V,$ then $mx_r\not\in VS_1.$  Otherwise,
$mx_r=nx_i$ for some $n\in V$ and some $i,$ and then $m>n$, a contradiction. To prove that  $\phi$  is
surjective, consider  a monomial $m$  in $S_{a+1}\setminus VS_1$.  Then $m=x_rn$ since $(x_0,\dots,
x_{r-1})^{a+1}\subset VS_1$. Obviously, $n\not\in V$. So $\phi$ is surjective. 
\end{proof}

\begin{prop}\label{ginpoints} Let $I$   be the ideal defining $s$ points, say  $P_1,\dots,
P_s$, of ${\bf P}^r$. Assume that there exists a coordinate system $x_0,x_1,\dots, x_r$ such that $I$
does not contain forms of degree $\leq s$ supported on $\leq s$ monomials.  Then  $\gin_\tau
I=\Seg_\tau(I)$ for all term orders $\tau$. In particular $\gin_{\lex} I=\Seg_{\lex}(I)$.
\end{prop}

\begin{proof} It is easy to see that the assumption implies that the Hilbert function of $S/I$ is the
expected one, namely  $\dim [S/I]_d= \min\{ s, {r+d \choose r}\}$ for all $d$.   Fix a term order
$\tau$.   For a given
$d\leq s$ consider  the set $M_d$  of the smallest (with respect to $\tau$) $\min\{ s, {r+d \choose
r}\}$  monomials of degree $d$. By assumption these monomials are a basis  of $S/I$ in degree $d$. It
follows immediately that  $\ini_\tau I_d=\Seg_{\tau}(I)_d$ for every $d\leq s$.   From Lemma \ref{1-dim}
we know that   $\ini_\tau I$ does not have generators in degree $\geq s$. Then $\ini_\tau I_d=\ini_\tau
I_sS_{d-s}$ for all $d\geq s$. On the other hand, it follows from Lemma \ref{segments} that
$\Seg_\tau(I)_d=\Seg_\tau(I)_sS_{d-s}$ for all $d\geq s$. We have seen already that
$\ini_\tau I_s=\Seg_\tau(I)_s$. Therefore  we may conclude that $\ini_\tau I_d=\Seg_\tau(I)_d$ also for
all
$d\geq s$. We have shown that $\ini_\tau I=\Seg_\tau(I)$. From this it follows that
$\gin_\tau I=\Seg_\tau(I)$ (see the construction/definition of gin given in \cite[Theorem 15.18]{eisenbud}.)
\end{proof}

We can now prove the main result of this section:

\begin{proof} [Theorem \ref{thm:points}] Let $P_1, \ldots, P_s$ be  generic  points in $\PP^r.$  Fix  a
coordinate system on
$\PP^r$ and let $(a_{i0},a_{i1},  \dots, a_{ir})$  be the coordinates of $P_i$. It is enough to show
that  the assumption of Proposition \ref{ginpoints} holds (in the given coordinates) for a generic
choice of the
$a_{ij}$. For any $d\leq s$ consider the $s\times {r+d \choose r}$ matrix $X_d$ whose rows are indexed 
by the points,  the  columns by the monomials of degree $d$ and whose $ij$-th entry is obtained by
evaluating the $j$-th monomial at the $i$-th point.  The assumption of Proposition \ref{ginpoints}  is
equivalent to the fact that any maximal minor of $X_d$ is non-zero for $d\leq s$. If we consider the
$a_{ij}$ as variables over some base field  then  every minor of $X_d$ is a non-zero polynomial in the
$a_{ij}$ since no cancellation can occur in the  expansion.   So these are finitely many non-trivial 
polynomial conditions on the coordinates of the points.   
\end{proof}

As we have already said, the genericity condition required in Theorem \ref{thm:points} implies that the 
Hilbert function of the ideal $I$ of $s$   points of $\PP^r$ is given is the expected one: 

$$\dim [S/I]_j=\min( s, {r+j \choose r}).$$

One may wonder whether it is enough to assume that the Hilbert function is generic to conclude that
$\gin_\tau I$ is $\Seg_\tau(I)$ for an ideal of points. The next example  answer this question. 

\begin{ex} (a) Consider the ideal $I$ of $7$ points of ${\bf P}^3$ with generic Hilbert function. The
ideal $I$ contains $3$ quadrics. If the $3$ quadrics have a common linear factor, then $\gin(I_2)$ is
$x_0(x_0,x_1,x_2)$ no matter what the term order is. So in particular, $\gin_{\revlex} I$ is not 
$\Seg_{\revlex}(I)$ in degree $2$. Explicitely,  one can take the $7$ points with coordinates 
$(0,0,0,1), (0,0,1,1), Ê(0,0,2,1), (0,1,0,1), (0,1,1,1), (0,2,0,1), (1,0,0,1)$. 

(b) Consider the $10$ points of ${\bf P}^3$ with coordinates $(a,b,c,1)$ where $a,b,c$ are non negative
integers with $a+b+c\leq 2$ and let $I$ the corresponding ideal.  One can check  with (and even without)
the help of a computer algebra system that the $10$ points have the generic Hilbert function and that
any generic projection to
${\bf P}^2$ gives $10$ points on a cubic. This, in turn, implies that $\gin_{\lex} I$ contains
$x_2^3$ while $\Seg_{\lex}(I)$ does not contain it. 
\end{ex} 

The next example shows that, even for Hilbert functions of generic  points in
$\PP^2$, the segment ideals are special among the Borel-fixed ideals.

\begin{ex}  Consider the ideal $I$ of $7$ generic points in $\PP^2$.  The Hilbert function of $S/I$ is
$(1,3,6,7,7,7,\dots)$. There are exactly $8$  Borel-fixed ideals with this Hilbert function, they are:

\[
\begin{array}{lll} (1) & (x^3, x^2y, x^2z, xy^3, xy^2z, xyz^3, xz^5, y^7), & \lex \\ (2) & (x^3, x^2y,
x^2z, xy^3, xy^2z, xyz^3, y^6), & (6, 2, 1)  \\ (3) &(x^3, x^2y, x^2z, xy^3, xy^2z, y^5), & (4, 2, 1)  \\
(4) & (x^3, x^2y, x^2z,  xy^3, y^4),  \\ (5) &(x^3, x^2y, xy^2, x^2z^2, xyz^3, xz^5, y^7), \\ (6) &(x^3,
x^2y, xy^2, x^2z^2, xyz^3, y^6),  \\ (7) &(x^3, x^2y, xy^2, x^2z^2, y^5), \\ (8) & (x^3, x^2y, xy^2,
y^4) & \revlex
\end{array}
\]

The ideals (1),(2),(3) and (8) are segments (with respect to the term  order or weight indicated on the
right) while the remaining four are  non-segments. Let us check, for instance, that (4) is not a segment.
Suppose, by  contradiction,  it is a segment with respect to  a term order $\tau.$ Then since
$x^2z$ is in and $xy^2$ is out, we have  $x^2z>_\tau xy^2$ and hence 
$xz>_\tau y^2.$  We deduce that $xy^2z>_\tau y^4.$ But  since $y^4$ is in then  also $xy^2z$ must be in
and this is a contradiction. Summing up, among the $8$  Borel-fixed ideals only (1),(2),(3) and (8) are
gins of $I.$
\end{ex}


\begin{thebibliography}{10}

\bibitem{bayer} D. Bayer, \emph{The Division Algorithm and the Hilbert Scheme}, Ph.D. thesis, Harvard University, 1982.

\bibitem{bayer-stillman} D.Bayer and M.Stillman, \emph{A criterion for detecting
  m-regularity}, Invent.Math. \textbf{87} (1987), 1--11.

\bibitem{bermejo-lejeune} I. Bermejo and M. Lejeune-Jalabert, \emph{Sur la 
complexit\'e du calcul des projections d'une courbe projective}, Comm. Alg.
\textbf{27} (1999), no. 7, 3211--3220.

\bibitem{bigatti} A. Bigatti, \emph{Upper bounds for the Betti numbers of a 
given Hilbert function}, Comm. Algebra \textbf{21} (1993), no. 7, 2317--2334.

\bibitem{bruns-herzog} W.Bruns and J.Herzog, \emph{Cohen Macaulay Rings}, Cambridge University Press,
Cabridge (1993).


\bibitem{bruns-herzog1}  W.Bruns and J.Herzog, 
\emph{ On the computation of $a$-invariants},  Manuscripta Math.\textbf{77} (1992), no.2-3, 201--213.

\bibitem{bruns-vetter} W.Bruns and U.Vetter, \emph{Determinantal Rings}, Lecture Notes in Mathematics,
\textbf{1327}, Springer-Verlag, New York, 1980.

\bibitem{CM} M. Chardin and G. Moreno-Soc\'ias, \emph{Regularity of lex-segment
ideals:some closed formulas and applications}, Proc.Amer.Math.Soc. \textbf{131} (2002),  no. 4, 1093--1102.

\bibitem{cocoa}A.Capani, G.Niesi, L.Robbiano, \emph{CoCoA, a system for doing Computations in
Commutative Algebra}, Available via anonymous ftp from: cocoa.dima.unige.it.

\bibitem{BudCasGor} N.Budur, M.Casanellas and E.Gorla \emph{ Hilbert functions of irreducible
arithmetically Gorenstein schemes}, to appear in J.Algebra.  
 

\bibitem{cox-little-oshea} D.Cox, J.Little, and D.O'Shea, \emph{Ideals, varieties, and algorithms},
Springer-Verlag, New York (1997).


\bibitem{eisenbud} D.Eisenbud, \emph{An Introduction to Commutative Algebra with a View Towards
Algebraic Geometry}, Springer-Verlag, New York (1995).

\bibitem{eliahou-kervaire} S.Eliahou and M.Kervaire, \emph{Minimal resolutions of some monomial ideals},
J.Algebra \textbf{129} (1990), 1--25
 
\bibitem{grayson-stillman} D.Grayson and M.Stillman, \emph{Macaulay 2 -- a system for computation in 
algebraic geometry and commutative algebra},
  http://www.math.uiuc.edu/Macaulay2, 1997.

\bibitem{green} M.Green, \emph{Generic initial ideals}, Six lectures on commutative algebra (Bellaterra,
1996), Progr.Math., \textbf{166}, Birkh\"auser, Basel, (1998), 119--186.

\bibitem{haiman} M.Haiman, \emph{Hilbert schemes, polygraphs and the {M}acdonald positivity
              conjecture}, J.Amer.Math.Soc., \textbf{14}, 2001, 941-1006.

\bibitem{hartshorne} R.Hartshorne, \emph{Algebraic {G}eometry}, Springer-Verlag, New York, 1977.


\bibitem{hulett} H. Hulett, \emph{Maximum Betti numbers of homogeneous ideals 
with a given Hilbert function}, Comm. Algebra \textbf{21} (1993), no. 7. 2335--2350.

\bibitem{kreuzer-robbiano}  M. Kreuzer and L. Robbiano, \emph{Computational commutative algebra}, Springer-Verlag, Berlin, 2000. 

\bibitem{marinari-ramella} M.G. Marinari and L. Ramella, \emph{Some properties
of Borel ideals}, J. Pure Appl. Algebra, \textbf{139} (1999), no. 1-3, 183--2000.

\bibitem{pardue} K. Pardue, \emph{Deformation classes of graded modules and
maximal Betti numbers}, Illinois J. Math., \textbf{40} (1996), no. 4, 564--585.
\end{thebibliography}
\end{document}